\documentclass[11pt,leqno]{article}
\usepackage{amsfonts,amssymb,latexsym,amsthm}
\usepackage{fullpage}
\usepackage{graphicx,epic,eepic}

\newtheorem{theorem}{Theorem}
\newtheorem{lemma}[theorem]{Lemma}
\theoremstyle{remark}
\newtheorem{remark}{Remark}

\newcommand{\N}{{\mathbb N}}
\newcommand{\R}{{\mathbb R}}
\newcommand{\C}{{\mathcal C}}
\newcommand{\ul}{\underline}
\newcommand{\ol}{\overline}

\title{A critical phenomenon for sublinear \\ elliptic equations in cone--like domains}
\author{
{\Large Vladimir Kondratiev}\\
Department of Mathematics\\
and Mechanics\\
Moscow State University\\
Moscow 119 899, Russia\\
{\tt kondrat@vnmok.math.msu.su}\\
\and
{\Large Vitali Liskevich}\\
School of Mathematics\\
University of Bristol\\
Bristol BS8 1TW\\
United Kingdom\\
{\tt v.liskevich@bristol.ac.uk}\\
\smallskip\and
{\Large Vitaly Moroz}\\
School of Mathematics\\
University of Bristol\\
Bristol BS8 1TW\\
United Kingdom\\
{\tt v.moroz@bristol.ac.uk}\\
\and
{\Large Zeev Sobol}\\
Department of Mathematics\\
Imperial College, London\\
London SW7 2AZ\\
United Kingdom\\
{\tt z.sobol@imperial.ac.uk}}

\date{}

\begin{document}
\maketitle

\begin{abstract}
We study positive supersolutions to an elliptic equation $(\ast)$
$-\Delta u=c|x|^{-s}u^p$, $\:p,s\in\R$, in cone--like domains in $\R^N$ ($N\ge 2$).
We prove that in the sublinear case $p<1$ there exists a critical exponent
$p_\ast<1$ such that equation $(\ast)$ has a positive supersolution if and only if $-\infty<p<p_\ast$.
The value of $p_\ast$ is determined explicitly by $s$ and the geometry of the cone.
\end{abstract}

\section{Introduction}

We study the existence and nonexistence of positive solutions and supersolutions
to the equation
\begin{equation}\label{*}
-\Delta u=\frac{c}{|x|^{s}}u^p\quad\mbox{in }\:\C_\Omega^\rho.
\end{equation}
Here $p\in\R$, $s\in\R$, $c>0$ and $\C_\Omega^\rho\subset\R^N$ ($N\ge 2$)
is an unbounded cone--like domain
$$\C_\Omega^\rho:=\{(r,\omega)\in\R^N:\:\omega\in\Omega,\:r>\rho\},$$
where $(r,\omega)$ are the polar coordinates in $\R^N$, $\rho>0$ and
$\Omega\subseteq S^{N-1}$ is a subdomain (a connected open subset) of the unit sphere $S^{N-1}$ in $\R^N$.
We say that $u\in H^1_{loc}(\C_\Omega^\rho)$ is a {\em supersolution} ({\em subsolution})
to equation (\ref{*}) if
$$\int_{\C_\Omega^\rho}\nabla u\cdot\nabla\varphi\:dx\ge(\le)
\int_{\C_\Omega^\rho}\frac{c}{|x|^s}u^p\varphi\:dx\quad\mbox{for all }\:
0\le\varphi\in C^\infty_0(\C_\Omega^\rho).$$
If $u$ is a sub and supersolution to (\ref{*}) then $u$
is said to be a {\em solution} to (\ref{*}).
By the weak Harnack inequality any nontrivial nonnegative supersolution to (\ref{*})
is positive in $\C_\Omega^\rho$.

We define {\em critical exponents} for equation (\ref{*}) by
$$p^\ast=p^\ast(\Omega,s)=
\inf\{p>1:\mbox{(\ref{*}) has a positive supersolution in $\C_\Omega^\rho$ for some $\rho>0$}\},$$
$$p_\ast=p_\ast(\Omega,s)=
\sup\{p<1:\mbox{(\ref{*}) has a positive supersolution in $\C_\Omega^\rho$ for some $\rho>0$}\}.$$
Set $p_\ast=-\infty$ if (\ref{*}) has no positive supersolution in $\C_\Omega^\rho$
for any $p<1$.

\begin{remark}
$(i)$
One can show that if $p<p_\ast$ or $p>p^\ast$ then (\ref{*}) has a positive solution in $\C_\Omega^\rho$
(see \cite{KLM} for the proof of the case $p>1$ and the proofs below for the case $p<1$).
The existence (or nonexistence) of positive (super) solutions
at the critical values $p_\ast$ and $p^\ast$ is a separate issue.

$(ii)$
Observe that in view of the scaling invariance of the Laplacian the critical
exponents $p_\ast$ and $p^\ast$ do not depend on $\rho>0$.

$(iii)$
We do not make any assumptions on the smoothness of the domain $\Omega\subseteq S^{N-1}$.
\end{remark}

Let $\lambda_1=\lambda_1(\Omega)\ge 0$ be the principal eigenvalue
of the Dirichlet Laplace--Beltrami operator $-\Delta_\omega$ on $\Omega$.
Let $\alpha_+\ge 0$ and $\alpha_-<0$ be the roots of the quadratic equation
$$\alpha(\alpha+N-2)=\lambda_1(\Omega).$$
In the {\em superlinear} case $p>1$ the value of the critical exponent is $p^\ast=1-\frac{2-s}{\alpha_-}$.
Moreover, if $s<2$ then (\ref{*}) has no positive supersolutions in the critical case $p=p^\ast$.
This has been proved by Bandle and Levine \cite{Bandle-Levine}, Bandle and Essen \cite{Bandle}
and Berestycki, Capuzzo--Dolcetta and Nirenberg \cite{BCN}
(see also \cite{KLM} for yet another proof of this result and
for  equations with measurable coefficients).

The {\em sublinear} case $p<1$ has been studied in
\cite{Brezis-Kamin,KLS}. From the result of Brezis and Kamin \cite{Brezis-Kamin}
it follows that for $p\in(0,1)$ equation (\ref{*}) has a bounded
positive solution in $\R^N$ if and only if $s>2$.
It has been proved in \cite{KLS} (amongst other things) that
for any $p\in(-\infty,1)$ equation
(\ref{*}) has a positive supersolution outside a ball in $\R^N$ if and only if $s>2$.

In this note, we discover a new critical phenomenon. Namely, we show that in sublinear case
equation (\ref{*}) exhibits a "non-trivial" critical exponent ($p_\ast>-\infty$) in cone-like domains.
The main result of the paper reads as follows.

\begin{theorem}\label{Main}
For $p\le 1$, the critical exponent for equation (\ref{*}) is
$p_\ast=\min\{1-\frac{2-s}{\alpha_+},1\}$.
If $p_\ast<1$ then (\ref{*}) has no positive supersolutions in the critical case $p=p_\ast$.
\end{theorem}

\begin{figure}[t]
\setlength{\unitlength}{0.0006in}
\begingroup\makeatletter\ifx\SetFigFont\undefined%
\gdef\SetFigFont#1#2#3#4#5{%
  \reset@font\fontsize{#1}{#2pt}%
  \fontfamily{#3}\fontseries{#4}\fontshape{#5}%
  \selectfont}%
\fi\endgroup%
{\renewcommand{\dashlinestretch}{30}

\begin{picture}(4824,4800)(0,-10)
\texture{88555555 55000000 555555 55000000 555555 55000000 555555 55000000 
           555555 55000000 555555 55000000 555555 55000000 555555 55000000 
           555555 55000000 555555 55000000 555555 55000000 555555 55000000 
           555555 55000000 555555 55000000 555555 55000000 555555 55000000 }
\shade\path(3012,3012)(4512,12)(12,12)(12,612)(3012,3012)

\path(12,1812)(4812,1812) 
\blacken\path(4692.000,1782.000)(4812.000,1812.000)(4692.000,1842.000)(4692.000,1782.000)
\path(2412,12)(2412,4762) 
\blacken\path(2442.000,4642.000)(2412.000,4762.000)(2382.000,4642.000)(2442.000,4642.000)

\dottedline{120}(12,3012)(4812,3012)
\dottedline{120}(3012,4762)(3012,12)

\thicklines
\path(12,612)(3012,3012)(4512,12)

\texture{aa777777 77aaaaaa aad5d5d5 d5aaaaaa aa777777 77aaaaaa aadd5ddd 5daaaaaa
         aa777777 77aaaaaa aad5d5d5 d5aaaaaa aa777777 77aaaaaa aa5ddd5d ddaaaaaa
         aa777777 77aaaaaa aad5d5d5 d5aaaaaa aa777777 77aaaaaa aadd5ddd 5daaaaaa
         aa777777 77aaaaaa aad5d5d5 d5aaaaaa aa777777 77aaaaaa aa5ddd5d ddaaaaaa }
\put(3012,3012){\shade\ellipse{80}{80}}
\put(4200,1812){\shade\ellipse{60}{60}}
\put(1510,1812){\shade\ellipse{60}{60}}
\put(3620,1812){\shade\ellipse{60}{60}}
\put(2412,2540){\shade\ellipse{60}{60}}

\put(4690,1662){\makebox(0,0)[lb]{\smash{{{\SetFigFont{8}{12.0}{\rmdefault}{\mddefault}{\updefault}$p$}}}}}
\put(2260,4650){\makebox(0,0)[lb]{\smash{{{\SetFigFont{8}{12.0}{\rmdefault}{\mddefault}{\updefault}$s$}}}}}
\put(2280,3070){\makebox(0,0)[lb]{\smash{{{\SetFigFont{8}{12.0}{\rmdefault}{\mddefault}{\updefault}$2$}}}}}
\put(2900,1650){\makebox(0,0)[lb]{\smash{{{\SetFigFont{8}{12.0}{\rmdefault}{\mddefault}{\updefault}$1$}}}}}
\put(3580,1990){\makebox(0,0)[lb]{\smash{{{\SetFigFont{8}{12.0}{\rmdefault}{\mddefault}{\updefault}$1\!-\!\frac{2}{\alpha_-}$}}}}}
\put(3980,1550){\makebox(0,0)[lb]{\smash{{{\SetFigFont{8}{12.0}{\rmdefault}{\mddefault}{\updefault}$\frac{N}{N-2}$}}}}}
\put(950,1990){\makebox(0,0)[lb]{\smash{{{\SetFigFont{8}{12.0}{\rmdefault}{\mddefault}{\updefault}$1\!-\!\frac{2}{\alpha_+}$}}}}}
\put(1850,2600){\makebox(0,0)[lb]{\smash{{{\SetFigFont{8}{12.0}{\rmdefault}{\mddefault}{\updefault}$2\!-\!\alpha_+$}}}}}
\put(512,412){\makebox(0,0)[lb]{\smash{{{\SetFigFont{8}{12.0}{\rmdefault}{\mddefault}{\itdefault}Nonexistence zone}}}}}
\put(312,3512){\makebox(0,0)[lb]{\smash{{{\SetFigFont{8}{12.0}{\rmdefault}{\mddefault}{\itdefault}Existence zone}}}}}

\end{picture}
\hfill
\begin{picture}(4824,3789)(0,-10)

\texture{88555555 55000000 555555 55000000 555555 55000000 555555 55000000 
           555555 55000000 555555 55000000 555555 55000000 555555 55000000 
           555555 55000000 555555 55000000 555555 55000000 555555 55000000 
           555555 55000000 555555 55000000 555555 55000000 555555 55000000 }
\shade\path(4812,4446)(3012,3012)(2148,4762)(4812,4762)(4812,4446)

\path(12,1812)(4812,1812)
\blacken\path(4692.000,1782.000)(4812.000,1812.000)(4692.000,1842.000)(4692.000,1782.000)
\path(2412,12)(2412,4762) 
\blacken\path(2442.000,4642.000)(2412.000,4762.000)(2382.000,4642.000)(2442.000,4642.000)

\dottedline{120}(12,3012)(4812,3012)
\dottedline{120}(3012,4762)(3012,12)

\dashline{120}(12,612)(3012,3012)
\dashline{120}(3012,3012)(4512,12)

\thicklines
\path(4812,4446)(3012,3012)(2148,4762)

\texture{aa777777 77aaaaaa aad5d5d5 d5aaaaaa aa777777 77aaaaaa aadd5ddd 5daaaaaa
         aa777777 77aaaaaa aad5d5d5 d5aaaaaa aa777777 77aaaaaa aa5ddd5d ddaaaaaa
         aa777777 77aaaaaa aad5d5d5 d5aaaaaa aa777777 77aaaaaa aadd5ddd 5daaaaaa
         aa777777 77aaaaaa aad5d5d5 d5aaaaaa aa777777 77aaaaaa aa5ddd5d ddaaaaaa }
\put(3012,3012){\shade\ellipse{80}{80}}
\put(4200,1812){\shade\ellipse{60}{60}}
\put(1510,1812){\shade\ellipse{60}{60}}
\put(3620,1812){\shade\ellipse{60}{60}}
\put(2412,4240){\shade\ellipse{60}{60}}

\put(4690,1662){\makebox(0,0)[lb]{\smash{{{\SetFigFont{8}{12.0}{\rmdefault}{\mddefault}{\updefault}$p$}}}}}
\put(2255,4640){\makebox(0,0)[lb]{\smash{{{\SetFigFont{8}{12.0}{\rmdefault}{\mddefault}{\updefault}$\sigma$}}}}}
\put(2280,3070){\makebox(0,0)[lb]{\smash{{{\SetFigFont{8}{12.0}{\rmdefault}{\mddefault}{\updefault}$2$}}}}}
\put(2900,1650){\makebox(0,0)[lb]{\smash{{{\SetFigFont{8}{12.0}{\rmdefault}{\mddefault}{\updefault}$1$}}}}}
\put(3550,1990){\makebox(0,0)[lb]{\smash{{{\SetFigFont{8}{12.0}{\rmdefault}{\mddefault}{\updefault}$1\!-\!\frac{2}{\alpha_-}$}}}}}
\put(3980,1550){\makebox(0,0)[lb]{\smash{{{\SetFigFont{8}{12.0}{\rmdefault}{\mddefault}{\updefault}$\frac{N}{N-2}$}}}}}
\put(990,1990){\makebox(0,0)[lb]{\smash{{{\SetFigFont{8}{12.0}{\rmdefault}{\mddefault}{\updefault}$1\!-\!\frac{2}{\alpha_+}$}}}}}
\put(1850,4200){\makebox(0,0)[lb]{\smash{{{\SetFigFont{8}{12.0}{\rmdefault}{\mddefault}{\updefault}$2\!-\!\alpha_-$}}}}}
\put(3050,4460){\makebox(0,0)[lb]{\smash{{{\SetFigFont{8}{12.0}{\rmdefault}{\mddefault}{\itdefault}Nonexistence zone}}}}}
\put(652,612){\makebox(0,0)[lb]{\smash{{{\SetFigFont{8}{12.0}{\rmdefault}{\mddefault}{\itdefault}Existence zone}}}}}

\end{picture}

}
\caption{Existence and nonexistence zones for equations (\ref{*}) (left) and (\ref{hat}) (right).
\label{fig}}
\end{figure}
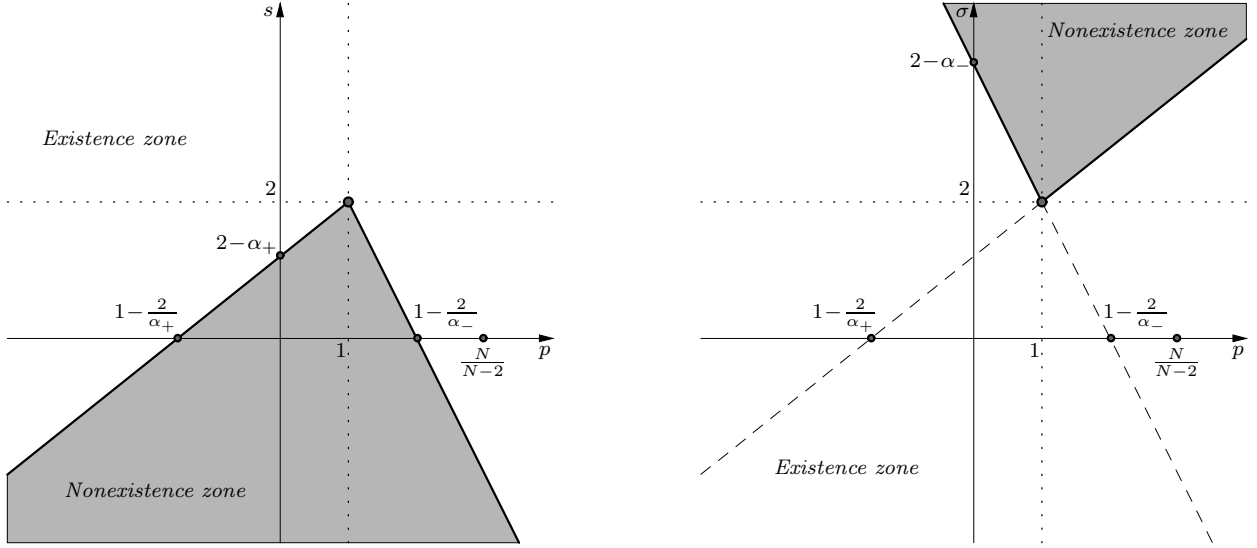

\begin{remark}
$(i)$ If $\alpha_+=0$ then we set $p_\ast=-\infty$.

$(ii)$
If $s>2$ then $p_\ast=p^\ast=1$ and (\ref{*}) has positive solutions for any $p\in\R$ \cite{Brezis-Kamin,KLS}.
If $s=2$ then $p_\ast=p^\ast=1$.
In this critical case (\ref{*}) becomes a linear equation with the
potential $c|x|^{-2}$, which has a positive (super) solution
if and only if $c\le \frac{(N-2)^2}{4}+\lambda_1(\Omega)$.

$(iii)$
Let $S_k=\{x\in S^{N-1}:x_1>0,\dots x_k>0\}$.
Then $\lambda_1(S_k)=k(k+N-2)$ and $\alpha_+(S_k)=k$, $\alpha_-(S_k)=2-N-k$.
Hence $p_\ast(S_k,s)=1-\frac{2-s}{k}$ and  $p^\ast(S_k,s)=1-\frac{2-s}{2-N-k}$.
In particular, in the case of the halfspace $S_1$ we have
$p_\ast(S_1,s)=s-1$ and $p^\ast(S_1,s)=\frac{N+1-s}{N-1}$.
\end{remark}

Applying the Kelvin transformation $y=y(x)=\frac{x}{|x|^2}$
we see that if $u$ is a positive solution to (\ref{*}) in $\C_\Omega^1$
then $\hat{u}(y)=|y|^{2-N}u(x(y))$ is a positive solution to
\begin{equation}\label{hat}
-\Delta\hat{u}=\frac{c}{|y|^{\sigma}}\hat{u}^p\quad\mbox{in }\:\widehat{\C}_\Omega^1,
\end{equation}
where $\sigma=(N+2)-p(N-2)-s$ and
$\widehat{\C}_\Omega^{\:1}:=\{(r,\omega)\in\R^N:\:\omega\in\Omega,\:0<r<1\}$.
We define the critical exponents $\widehat{p}^\ast=\widehat{p}^\ast(\Omega,s)$ and $\widehat{p}_\ast=\widehat{p}_\ast(\Omega,s)$
for equation (\ref{hat}) similarly to $p^\ast(\Omega,s)$ and $p_\ast(\Omega,s)$.
In the superlinear case $p>1$, Bandle and Essen \cite{Bandle} proved that if
$\sigma>2$ then
$\widehat{p}^\ast=1-\frac{2-\sigma}{\alpha_+}$ and (\ref{hat})
has no positive supersolutions when $p=\widehat{p}^\ast(\Omega)$.
In the sublinear case $p<1$ by an easy computation we derive from Theorem \ref{Main} the following result.

\begin{theorem}\label{Kelvin}
For $p\le 1$, the critical exponent for equation (\ref{hat}) is
$\widehat{p}_\ast=\min\{1-\frac{2-\sigma}{\alpha_-},1\}$.
If $\widehat{p}_\ast<1$ then (\ref{hat}) has no positive supersolutions
in the critical case $p=\widehat{p}_\ast$.
\end{theorem}

In the remaining part of the paper we prove Theorem \ref{Main}.

\section{Proof of Theorem \ref{Main}}

\paragraph{Existence.}
In the polar coordinates equation (\ref{*}) reads as follows
\begin{equation}\label{polar}
-u_{rr}-\frac{N\!-\!1}{r}u_r-\frac{1}{r^2}\Delta_\omega u = \frac{c}{r^{s}}u^p
\quad\mbox{in }\:\C_\Omega^1.
\end{equation}
Let $s\le 2, \, p< 1 - \frac{2-s}{\alpha_+}$.
Let $0<\psi\in H^1_{loc}(\Omega)$ be a positive solution to the equation
\begin{equation}\label{psi}
-\Delta_\omega\psi-\alpha(\alpha+N-2)\psi=\psi^{p}\quad\mbox{in }\:\Omega,
\end{equation}
where $\alpha:=\frac{2-s}{1-p}$.
Then  it is readily seen that
$u:=c^\frac{1}{1-p}r^{\alpha}\psi\in H^1_{loc}(\C_\Omega^1)$
is a positive solution to $(\ref{polar})$ in $\C_\Omega^1$.
Thus the problem reduces to the existence of
positive solutions to (\ref{psi}).

Note that
$0<\alpha(\alpha+N-2)<\lambda_1(\Omega)$.
Hence the operator $-\Delta_\omega-\alpha(\alpha+N-2)$
is coercive on $H^1_0(\Omega)$ and satisfies the maximum principle.
We consider separately the cases $p\in[0,1)$ and $p<0$.

\paragraph{\it Case $p\in[0,1)$.} 
Let $\phi_1>0$ be the principal Dirichlet eigenfunction of $-\Delta_\omega$ on $\Omega$.
Let $\ol\phi>0$ be the unique solution to the problem
$$-\Delta_\omega\phi-\alpha(\alpha+N-2)\phi=1,\qquad\phi\in H^1_0(\Omega).$$
Observe that $\phi_1, \ol\phi \in L^\infty$.

Hence
$\tau\ol\phi$ is a supersolution to (\ref{psi}) for a large $\tau>0$,
and $\epsilon\phi_1$ is a subsolution to (\ref{psi}) for a small $\epsilon>0$.
Thus by the sub and supersolutions argument equation (\ref{psi})
has a solution $\psi\in H^1_0(\Omega)$ such that $\epsilon\phi_1<\psi\le\tau\ol\phi$.

\paragraph{\it Case $p<0$.}
Consider the problem
\begin{equation}\label{psi-1}
-\Delta_\omega\phi-\alpha(\alpha+N-2)(\phi+1)=
(\phi+1)^{p},
\qquad\phi\in H^1_0(\Omega).
\end{equation}
Let $\ol\phi>0$ be the unique solution to the problem
$$-\Delta\phi-\alpha(\alpha+N-2)(\phi+1)=1,\qquad\phi\in H^1_0(\Omega).$$
It is clear that $\ol\phi$ is a supersolution to (\ref{psi-1}) and
$\ul\phi\equiv 0$ is a subsolution to (\ref{psi-1}).
We conclude that (\ref{psi-1}) has a positive solution $\phi\in H^1_0(\Omega)$
such that $0<\phi\le\ol\phi$.
Then $\psi:=\phi+1\in H^1_{loc}(\Omega)$ is a positive solution to (\ref{psi}).
This completes the proof of the existence part of Theorem \ref{Main}.

\paragraph{Nonexistence.}
In what follows we set $\delta:=1$ if $p<0$ and $\delta:=0$ if
$p\in[0,1)$. Let $G\subset\R^N$ be a domain, $0\not\in G$. Observe
that equation (\ref{*}) has a positive supersolution in $G$ if and
only if the equation
\begin{equation}\label{**}
-\Delta w=\frac{c}{|x|^{s}}(w+\delta)^p\quad\mbox{in }\:G
\end{equation}
has a positive supersolution.
Indeed, if $u>0$ is a supersolution to (\ref{*}) in $G$ then $u$
is a supersolution to (\ref{**}). If $w>0$ is a supersolution to (\ref{**})
then $u=w+\delta$ is a supersolution to (\ref{*}).
The main argument of the proof nonexistence rests upon the following two lemmas.

The next lemma is an adaptation a comparison principle
by Ambrosetti, Brezis and Cerami \cite[Lemma 3.3]{ABC}.
\begin{lemma}\label{Brezis}
Let $G\subset\R^N$ be a bounded domain, $0\not\in G$.
Let $0\le \ul{w}\in H^1_0(G)$ be a subsolution
and $0\le \ol{w}\in H^1_{loc}(G)$ a supersolution to (\ref{**}).
Then $\ul{w}\le\ol{w}$ in $G$.
\end{lemma}

\proof
In \cite[Lemma 3.3]{ABC} the result was proved for a smooth bounded domain $G$
and $\ul{w},\ol{w}\in H^1_0(G)$ (and more general nonlinearities).
The proof given in \cite{ABC} carries over literally to the case of
an arbitrary bounded domain $G$ and $\ul{w},\ol{w}\in H^1_0(G)$,
or a smooth bounded domain $G$,
$\ul{w}\in H^1_0(G)$  and $0\le\ol{w}\in H^1(G)$.
Thus we only need to extend the lemma to an arbitrary bounded domain $G$ and $\ol{w}\in H^1_{loc}(G)$.

Let $\ol w\in H^1_{loc}(G)$ be a supersolution to (\ref{**}) in $G$.
Let $(G_n)_{n\in\N}$ be an exhaustion of $G$, that is a sequence of bounded smooth domains
such that $\ol G_n\subset G_{n+1}\subset G$ and $\cup_{n\in\N}G_n=G$.
Analogously to the argument given above in the existence part of the proof,
one can readily see that, for each $n\in\N$,
there exists a solution $0<w_n\in H^1_0(G_n)$ to (\ref{**})
(e.g., by constructing appropriate sub and supersolutions).
Moreover, $w_n\le w_{n+1}$.
Observe that $w_n\le\ol w$ in $G_n$ by \cite[Lemma 3.3]{ABC}.

We claim that $\sup\|\nabla w_n\|_{L^2}<\infty$.
This is clear for $p<0$, since $(w+1)^p\le 1$.
For $p\in[0,1)$,
using $w_n$ as a test function in (\ref{**}), we have
$$\int_G|\nabla w_n|^2dx=\int_G\frac{c}{|x|^{s}}w_n^{p+1}\,dx\le
c_1\left(\int_G|\nabla w_n|^2 dx\right)^{(p+1)/2},$$
which implies the claim.
It follows that $w_n$ converges pointwise in $G$, strongly in $L^2(G)$ and weakly in $H^1_0(G)$
to a positive $w_\ast\in H^1_0(G)$.
Clearly $w_\ast>0$ is a solution to (\ref{**}) in $G$ and $0<w_\ast\le\ol w$ in $G$.

Now let $0\le \ul w\in H^1_0(G)$ be a subsolution to (\ref{**}) in $G$.
By \cite[Lemma 3.3]{ABC} we conclude that $\ul w\le w_\ast$ in $G$.
\qed

Next, consider the initial value problem
\begin{equation}\label{ODE}
-v_{rr}-\frac{N-1}{r}v_r+\frac{\lambda_1}{r^2}v=\frac{c}{r^s}v^p\quad\mbox{for }\:r>1;
\qquad v(1)=\delta,\quad v_r(1)=K;
\end{equation}
where $p<1$, $s\in\R$, $c>0$, $K>1$ and $\delta$ as above.
Let $(1,R)$, $R=R(\delta,K)\le\infty$,
be the maximal right interval of existence of the solution
$v$ to (\ref{ODE}) in the region $\{(r,v)\in(1,+\infty)\times(\delta,+\infty)\}$.

\begin{lemma}\label{L-ODE}
Let $s<2$ and $p\in [1-\frac{2-s}{\alpha_+},1)$.
Then for any interval $[r_\ast,r^\ast]\subset(1,+\infty)$ there exists $K_0>1$ such that
\begin{itemize}
\item[$i)$]
for all $K>K_0$ one has $r^\ast<R<+\infty$ and $v(r)\to\delta$ as $r\nearrow R$;
\item[$ii)$]
for any $M>\delta$ there exists $K>K_0$ such that $\min_{[r_\ast,r^\ast]}v\ge M$.
\end{itemize}
\end{lemma}



\proof
Set $\alpha:=\alpha_+$, $v:=wr^{\alpha}$, $t=r^{2-N-2\alpha}$.
Then $w$ solves the following problem
$$w_{tt}+c_1t^{-\sigma}w^p=0\quad\mbox{for }\:t\in(T,1);\qquad w(1)=\delta,\quad w_t(1)=-L,$$
where $\sigma=\frac{2N-2+\alpha(p+3)-s}{N-2+2\alpha}\ge 2$, $c_1>0$, $0\le T=R^{2-N-2\alpha}<1$
and $L=\frac{K-\alpha\delta}{N-2+2\alpha}\to\infty$ as $K\to\infty$. Choose $K_0$ such that $L>\delta$.
Observe that $w(t)$ is concave, hence
$$\delta<w(t)\le w(1)-w_t(1)(1-t)\le \delta+L\quad\mbox{for }\:t\in(T,1).$$
To see that $T>0$ let $\tilde w:=w$ for $p<0$, otherwise let $\tilde w:=w^{1-p}$.
Then $\tilde w$ satisfies the inequality
$$\tilde w_{tt} + c_2 t^{-2}\tilde w^{q}\le 0\quad\mbox{for }\:t\in(T,1),$$
with $c_2>0$ and $q:=\min\{p,0\}$.
Integrating $\tilde w_{tt}$ twice one can easily see that such inequality
has no positive solutions in any neighborhood of zero.
Thus we conclude that $T>0$, hence $w(t)\to\delta$ as $t\searrow T$.
In particular, $w(t)$ attains its maximum on $(T,1)$.

Let $T_0\in(T,1)$ be such that $w_t(T_0)=-\frac{L-\delta}{2}$.
Since $\delta\le w(t)\le\delta + L$ for $t\in(T_0,1)$, it follows that
$$\frac{L+\delta}{2}=w_t(T_0)-w_t(1)=-\int_{T_0}^{1}w_{tt}d\tau=
c_1\int_{T_0}^{1}\frac{w^p}{\tau^\sigma}d\tau\le c_3\left(\frac{1}{T_0^{\sigma-1}}-1\right)
\quad\mbox{for }\:t\in(T_0,1).$$
Hence $T_0\to 0$ as $L\to+\infty$.
Therefore for any given $t^\ast<1$ there exists $L_0>1$ such that for any $L>L_0$ one has $0<T<T_0<t^\ast$.
Thus, $(i)$ follows with $r^\ast=(t^\ast)^{\frac{1}{N-2+2\alpha}}$.

Observe now that for any $L>L_0$ we have
$$-\frac{L-\delta}2\ge w_t(t)\ge -L \quad\mbox{for }\:t\in(t^\ast,1),$$
since $w$ is concave. Hence for any $t\in(t^\ast,1)$ we obtain
$$w(t)=w(1)-\int_{t}^{1}w_t\,d\tau\ge\delta+(1-t)\frac{L-\delta}{2}\to\infty\quad
\mbox{as }\:L\to\infty.$$
Thus $(ii)$ follows.
\qed

\paragraph{\it Nonexistence -- completed.}

Let $p\in[1-\frac{2-s}{\alpha_+},1)$.
Fix a compact $K\subset \C^1_\Omega$ and $M>1$. There exists an interval
$[r_\ast,r^\ast]\subset(1,+\infty)$ such that $K\subset \C_\Omega^{(r_\ast,r^\ast)}$,
where $\C_\Omega^{(r_1,r_2)}$ denotes the set $\{x\in\C_\Omega^1\,|\,r_1\le |x|\le r_2\}$.
Then by  Lemma \ref{L-ODE}
there exists
$v:(1,R)\to(\delta,+\infty)$ solving (\ref{ODE}) such that $R>r^\ast$ and $\inf_{[r_\ast,r^\ast]}v\ge M+\delta$.

Let $\phi_1>0$ be the principal Dirichlet eigenvalue of $-\Delta_\omega$ on $\Omega$
with $\|\phi_1\|_\infty =1$.
Set $w_M:=(v-\delta)\phi_1$. Then $0<w_M\in H^1_0(\C_\Omega^{(1,R)})$,
and direct computation shows that $w_M$ is a subsolution to (\ref{**}) in $\C_\Omega^{(1,R)}$.
Now assume that $w>0$ is a supersolution to (\ref{**}) in $\C_\Omega^1$.
By Lemma \ref{Brezis} it follows that that $w\ge w_M$ in $\C_\Omega^{(1,R)}$.
By the weak Harnack inequality we have
$$\inf_K w\ge c_H\int_K w\,dx\ge c_H\int_K w_M\,dx\ge c_2M.$$
Since $M$ was arbitrary, we conclude that $w\equiv+\infty$ in $K$.
\qed
\medskip

\begin{small}
{\bf Acknowledgement.}
The support of the MFI Oberwolfach, Nuffield Foundation
and Institute of Advanced Studies of Bristol University is gratefully acknowledged.
\end{small}

\begin{small}

\end{small}

\end{document}